\documentclass[a4paper,12pt]{article}

\usepackage{makeidx}

\usepackage{latexsym}
\usepackage{amscd}
\usepackage{amsmath}
\usepackage{amssymb}

\usepackage{amsthm}
\usepackage{float}

\usepackage{psfrag}

\theoremstyle{plain}
\newtheorem{thm}{Theorem}[section]
\newtheorem{theorem}[thm]{Theorem}

\newtheorem{proposition}[thm]{Proposition}

\theoremstyle{definition}
\newtheorem{definition}[thm]{Definition}
\newtheorem{remark}[thm]{Remark}

\newcommand{\suchthat}{\; | \;}
\renewcommand{\Re}{\ensuremath{\mathfrak{Re}}}
\renewcommand{\Im}{\ensuremath{\mathfrak{Im}}}
\newcommand{\Ker}{\mathrm{Ker}}

\newcommand{\Ob}{\mathfrak{Ob}}
\renewcommand{\hat}{\widehat}
\newcommand{\Hess}{\mathrm{Hess}}
\newcommand{\coh}{\ensuremath{\mathrm{coh}}}
\newcommand{\Ext}{\mathop{\mathrm{Ext}}\nolimits}
\newcommand{\Hom}{\mathop{\mathrm{Hom}}\nolimits}
\newcommand{\RHom}{\mathop{\bR\mathrm{Hom}}\nolimits}
\newcommand{\RGamma}{\mathop{\bR\Gamma}\nolimits}

\newcommand{\End}{\mathop{\mathrm{End}}\nolimits}
\newcommand{\Spec}{\mathrm{Spec}}

\newcommand{\M}{\overline{\scM}}

\newcommand{\Fuk}{\mathfrak{Fuk}}
\newcommand{\m}{\mathfrak{m}}

\newcommand{\Lag}{\mathcal{L}ag}

\renewcommand{\det}{\mathop{\mathrm{det}}\nolimits}

\newcommand{\bC}{\ensuremath{\mathbb{C}}}

\newcommand{\bP}{\ensuremath{\mathbb{P}}}

\newcommand{\bR}{\ensuremath{\mathbb{R}}}

\newcommand{\bZ}{\ensuremath{\mathbb{Z}}}

\newcommand{\scA}{\ensuremath{\mathcal{A}}}

\newcommand{\scC}{\ensuremath{\mathcal{C}}}

\newcommand{\scE}{\ensuremath{\mathcal{E}}}
\newcommand{\scF}{\ensuremath{\mathcal{F}}}

\newcommand{\scI}{\ensuremath{\mathcal{I}}}

\newcommand{\scM}{\ensuremath{\mathcal{M}}}

\newcommand{\scO}{\ensuremath{\mathcal{O}}}

\newcommand{\scS}{\ensuremath{\mathcal{S}}}

\newcommand{\scU}{\ensuremath{\mathcal{U}}}
\newcommand{\scV}{\ensuremath{\mathcal{V}}}

\newcommand{\Res}{\ensuremath{\mathop{\mathrm{Res}}}}
\newcommand{\Gr}{\ensuremath{\mathrm{Gr}}}
\newcommand{\id}{\ensuremath{\mathop{\mathrm{id}}}}
\newcommand{\Sym}{\ensuremath{\mathop{\mathrm Sym}\nolimits}}
\newcommand{\econe}{\ensuremath{\Lphi \scO_{\bP^2}(-1)}}
\newcommand{\ectwo}{\ensuremath{\Lphi \Omega_{\bP^2}(1)}}
\newcommand{\ecthree}{\ensuremath{\scO_\dP}}
\newcommand{\eeone}{\ensuremath{\scE_1}}
\newcommand{\eetwo}{\ensuremath{\scE_2}}
\newcommand{\eethree}{\ensuremath{\scE_3}}
\newcommand{\eefour}{\ensuremath{\scE_4}}
\newcommand{\eefive}{\ensuremath{\scE_5}}
\newcommand{\eesix}{\ensuremath{\scE_6}}

\newcommand{\Lphi}{\phi^*}
\newcommand{\Rphi}{\bR \phi_*}
\newcommand{\rt}{{\small \mathrm{right}}}
\newcommand{\lt}{{\small \mathrm{left}}}
\newcommand{\modm}{\mathrm{mod}\mbox{-}}
\newcommand{\dP}{Y}
\newcommand{\DbdP}{\ensuremath{{D^b \mathrm{coh}(Y)}}}
\newcommand{\flata}{\ensuremath{\alpha}}

\newcommand{\canonicala}{\ensuremath{i}}
\newcommand{\canonicalb}{\ensuremath{j}}
\newcommand{\FukW}{\Fuk^{\rightarrow} W}
\newcommand{\hFukW}{{\hom}_{\FukW}}
\newcommand{\HFukW}{\Hom_{\FukW}}
\newcommand{\vspan}{\ensuremath{\mathop{\mathrm{span}}}}
\title{Homological Mirror Symmetry for Toric del Pezzo Surfaces}
\author{Kazushi Ueda}
\date{}
\begin{document}

\maketitle

\begin{abstract}

We prove the homological mirror conjecture for toric del Pezzo surfaces.
In this case, the mirror object is a regular function
on an algebraic torus $(\bC^\times)^2$.
We show that the derived Fukaya category of this mirror
coincides with the derived category of coherent sheaves
on the original manifold.

\end{abstract}

\section{Introduction}

Mirror symmetry
started as a mysterious relationship
between complex geometry of a Calabi-Yau 3-fold
and symplectic geometry of another Calabi-Yau 3-fold
called the {\em mirror manifold}.
%See e.g. \cite{Cox-Katz} for a general background
%on mirror symmetry.
In 1994, Kontsevich \cite{Kontsevich_HAMS} proposed
a program to understand various mirror phenomena
as a consequence of the following
{\em homological mirror conjecture}:
Calabi-Yau manifolds always come in pairs
in such a way that
the derived category of coherent sheaves
on one manifold is equivalent
as a triangulated category
to the derived Fukaya category of the other.

Although mirror symmetry was first discovered
for Calabi-Yau manifolds,
there are also variants of these phenomena
for other classes of manifolds.
One such example is the Givental's theorem
\cite{Givental_HGMS}
giving integral representations of
$J$-functions
of toric Fano manifolds.
See also \cite{Batyrev_QCRTM}, \cite{Witten_PN2T2D},
\cite{Eguchi-Hori-Xiong} and \cite{Hori-Vafa}.

We can also formulate the homological mirror conjecture
for toric Fano manifolds.
Let $X$ be a toric Fano manifold of dimension $n$.
Then the mirror partner
is
a regular function $W$
on
an algebraic torus $(\bC^\times)^n$
of the same dimension
equipped with a symplectic structure
(along with an additional data
called a {\em grading}).
Here, the function $W$ is a Newton polynomial
for the convex hull of the generators
of 1-dimensional cone of the fan of $X$.
Coefficients of this polynomial do not matter
as long as they are chosen general enough.
By Kouchnirenko \cite{Kouchnirenko},
$W$ has exactly $\dim H^*(X,\bC)$ critical points.
Take a regular value $t$ of $W$ and
a distinguished basis of
vanishing cycles in $W^{-1}(t)$.
We also have to choose
a {\em grading} and
a {\em spin structure}
on each of these vanishing cycles.
The {\em directed Fukaya category}
$\FukW$
of $W$
(along with the choice of gradings
and spin structures)
is an $A_\infty$-category
whose objects are vanishing cycles
and whose morphisms are Floer complexes.
Roughly speaking, the Floer complex
between two vanishing cycles
are the vector space
spanned by intersection points
between them,
and the compositions of morphisms
are given by ``counting polygons.''
By Seidel \cite{Seidel_VC},
the derived category
$D^b \FukW$ of
$\FukW$
is independent of the choice
of a distinguished basis
of vanishing cycles.

The following is the main result of this paper:

\begin{theorem} \label{th:hms}
The derived category of coherent sheaves
on a toric del Pezzo surface $X$
is equivalent as a triangulated category
to the derived category of
the directed Fukaya category
of the mirror $W$ of $X$;
$$
 D^b \coh(X) \cong D^b \FukW.
$$
\end{theorem}

Our proof is based on
an explicit computation.
The above Theorem \ref{th:hms}
extends the work of
Seidel \cite{Seidel_VC2}
and 
Auroux, Katzarkov and Orlov
\cite{Auroux-Katzarkov-Orlov},
where the cases of
$\bP^2$ and $\bP^1 \times \bP^1$,
and
$\bP^2$ blown-up at one point
is treated.
See also the paper
by Hori, Iqbal and Vafa
\cite{Hori-Iqbal-Vafa},
where mirror symmetry
for Fano manifolds
are discussed
from a physics point of view.

{\bf Acknowledgements}:
We thank O. Fujino, K. Fukaya,
H. Iritani, T. Kawai,
K. Saito and A. Takahashi
for valuable discussions and comments.
The author is supported by
JSPS Fellowships for Young Scientists
No.15-5561.

\section{Derived category of coherent sheaves}

\begin{figure}[H]
\centering
\psfrag{v1}{$v_1$}
\psfrag{v2}{$v_2$}
\psfrag{v3}{$v_3$}
\psfrag{v4}{$v_4$}
\psfrag{v5}{$v_5$}
\psfrag{v6}{$v_6$}
\psfrag{vequalssomething}{$\begin{array}{rcl}
 v_1 & = & (   1   ,   0   ), \\
 v_2 & = & (   0   ,   1   ), \\
 v_3 & = & (   -1  ,  -1   ), \\
 v_4 & = & (   -1  ,   0   ), \\
 v_5 & = & (   0   ,  -1   ), \\
 v_6 & = & (   1   ,   1   ).
	     \end{array}$}
\includegraphics{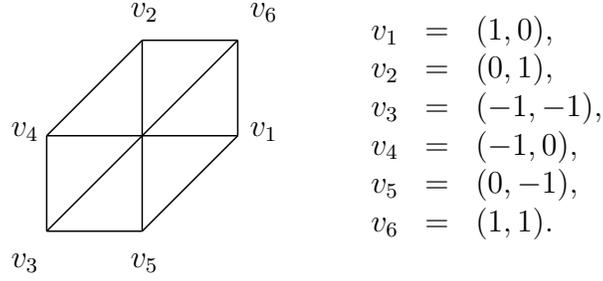}
\caption{the toric data of $\dP$}
\label{fg:dp3}
\end{figure}

We describe the structure
of the bounded derived category $\DbdP$
of coherent sheaves
on $\dP$ in this section,
where $\dP$ is the projective plane
blown-up
at three points $p_1$, $p_2$ and $p_3$
in general position.
Let $\phi : \dP \rightarrow \bP^2$ be this blow-up
and 
$E_1$, $E_2$ and $E_3$
be the exceptional divisors
corresponding to
$p_1$, $p_2$ and $p_3$ respectively.
It is a toric surface and
the generators of one-dimensional cones
of its fan is drawn
in Figure \ref{fg:dp3}.
%Then we describe the relation between
%$\DbdP$ and the derived category of
%coherent sheaves
%on other toric del Pezzo surfaces.

\begin{definition} \label{def:exceptional_collection}
\begin{enumerate}
\item 
 An object $\scE$ in a triangulated category
 is exceptional if
  $$
   \Ext^i(\scE,\scE)=\left\{\begin{array}{cl}
      \bC & \mbox{if $i=0$,}\\
       0  & \mbox{otherwise.}\\
   \end{array}\right.
  $$
\item
 An ordered set of objects $(\scE_i)_{i=1}^N$
 in a triangulated category
 is an exceptional collection
 if each $\scE_i$ is exceptional and
 $\Ext^k(\scE_i,\scE_j)=0$ for any $i>j$
 and for any $k$.
\end{enumerate}

\end{definition}

By combining theorems of Beilinson \cite{Beilinson} and
Orlov \cite{Orlov},
we have the following generators
of $\DbdP$:

\begin{theorem}
Let
\begin{eqnarray*}
\scC & = &
( \scE_1, \scE_2, \scE_3, \scE_4, \scE_5, \scE_6 )
\end{eqnarray*}
where
\begin{equation} \label{eq:ec}
 \begin{array}{rclrclrcl}
  \scE_1 & = & \scO_{E_1}(-1)[-1], &
  \scE_2 & = & \scO_{E_2}(-1)[-1], &
  \scE_3 & = & \scO_{E_3}(-1)[-1], \\
  \scE_4 & = & \econe, & 
  \scE_5 & = & \ectwo, &
  \scE_6 & = & \ecthree. \\
 \end{array}
\end{equation}
Then $\scC$ is an exceptional collection
generating $\DbdP$.
\end{theorem}
Here, $\scO_{E_i}(-1)$ is the sheaf supported on
$E_i$ and
is isomorphic to the tautological sheaf $\scO_{\bP^1}(-1)$
on $E_i \cong \bP^1$,
$[\bullet]$ is the shift operator in the derived category,
$\Lphi$ is the derived pull-back,
$\scO_{\bP^2}(-1)$ is the tautological sheaf
on $\bP^2$,
and $\Omega_{\bP^2}(1)$ is the cotangent sheaf of $\bP^2$
tensored with the hyperplane sheaf
$\scO_{\bP^2}(1)$.
%Note that $\scO_{\bP^2}(-1)$ is the exterior power
%$\Omega^2_{\bP^2}(2)$ of
%$\Omega_{\bP^2}(1)$.
%$N$ is six in our case.
%Since $\scC$ generates $\DbdP$,
%its categorical structure
%is wholly determined
%by the structure of morphisms
%between objects
%in $\scC$.

\begin{proposition}
All the non-zero $\Ext$-groups
within the exceptional collection $\scC$ are
% \begin{eqnarray*}
%  \Hom(\scO_{E_i}(-1), \econe) 
%   & = & \bC p_i, \\
%  \Hom(\scO_{E_i}(-1), \ectwo)
%   & = & \Ker(V^\vee \rightarrow (\bC p_i)^\vee), \\
%  \Hom(\scO_{E_i}(-1), \ecthree)
%   & = & \bC, \\
%  \Hom(\econe, \ectwo)
%   & = & \Lambda^2 V^\vee, \\
%  \Hom(\econe, \ecthree)
%   & = & V^\vee, \\
%  \Hom(\ectwo, \ecthree)
%   & = & V,
% \end{eqnarray*}
 \begin{eqnarray*}
  \Hom(\scE_i, \scE_4)
   & = & \bC p_i, \\
  \Hom(\scE_i, \scE_5)
   & = & \Ker(V^\vee \rightarrow (\bC p_i)^\vee), \\
  \Hom(\scE_i, \scE_6)
   & = & \bC, \\
  \Hom(\scE_4, \scE_5)
   & = & \Lambda^2 V^\vee, \\
  \Hom(\scE_4, \scE_6)
   & = & V^\vee, \\
  \Hom(\scE_5, \scE_6)
   & = & V,
 \end{eqnarray*}
where $i=1,2,3$,
$V \cong \bC^3$ is the three-dimensional vector space such that
$\bP^2 = \bP(V)$,
the check denotes the dual vector space,
$\bC p_i \subset V$ is the one-dimensional subspace
corresponding to $p_i \in \bP(V)$,
and the map $V^\vee \rightarrow (\bC p_i)^\vee$
is the dual of the inclusion
$\bC p_i \hookrightarrow V$.
Compositions of morphisms
are given by
$$
\begin{array}{ccc}
 \Hom(\scE_i, \eefour) \times \Hom(\eefour, \eefive)
   & \longrightarrow & \Hom(\scE_i, \eefive)\\
 \rotatebox{90}{$\in$} & & \rotatebox{90}{$\in$} \\
 (v, \omega ) & \longmapsto & -\iota_v \omega
\end{array},
$$
$$
\begin{array}{ccc}
 \Hom(\scE_i, \eefour) \times \Hom(\eefour, \eesix)
   & \longrightarrow & \Hom(\scE_i, \eesix)\\
 \rotatebox{90}{$\in$} & & \rotatebox{90}{$\in$} \\
 (v, \omega ) & \longmapsto & \omega(v)
\end{array},
$$
$$
\begin{array}{ccc}
 \Hom(\scE_i, \eefive) \times \Hom(\eefive, \eesix)
   & \longrightarrow & \Hom(\scE_i, \eesix)\\
 \rotatebox{90}{$\in$} & & \rotatebox{90}{$\in$} \\
 (\omega, v ) & \longmapsto & \omega(v)
\end{array},
$$
$$
\begin{array}{ccc}
 \Hom(\eefour, \eefive) \times \Hom(\eefive, \eesix)
   & \longrightarrow & \Hom(\eefour, \eesix)\\
 \rotatebox{90}{$\in$} & & \rotatebox{90}{$\in$} \\
 (\omega, v ) & \longmapsto & \iota_v \omega
\end{array}
$$
where $\iota_v : \wedge^2 V^\vee \rightarrow V^\vee$
for $v \in V$ is the interior product.
\end{proposition}

\proof
All the computations reduce to those on $\bP^2$
in the following way:
Let $\scE$ and $\scF$ be vector bundles
on $\bP^2$. Then
\begin{eqnarray*}
 \RHom(\Lphi \scE, \Lphi \scF)
 & = & \RHom(\scE, \Rphi \Lphi \scF) \\
 & = & \RHom(\scE, \scF \otimes \Rphi \scO_\dP) \\
 & = & \RHom(\scE, \scF),
\end{eqnarray*}
\begin{eqnarray*}
 \RHom(\Lphi \scE, \scO_{E_i}(-1))
 & = & \RHom(\scE, \Rphi \scO_{E_i}(-1)) \\
 & = & 0,
\end{eqnarray*}
and
\begin{eqnarray*}
\RHom(\scO_{E_i}(-1)[-1], \Lphi \scE)
 & = & \RHom( \{ \scO_\dP \rightarrow \scO_\dP(E_i) \},
   \Lphi \scE ) \\
 & = & \RGamma( \{ \scO_\dP(-E_i) \rightarrow \scO_\dP \} \otimes
   \Lphi \scE ) \\
 & = & \RGamma( \{ \Rphi \scO_\dP(-E_i) \rightarrow
   \Rphi \scO_\dP \} \otimes \scE ) \\
 & = & \RGamma( \{ \scI_{p_i} \rightarrow \scO_{\bP^2} \}
   \otimes \scE ) \\
 & = & \RGamma( \scO_{p_i} \otimes \scE ) \\
 & = & \scE|_{p_i},
\end{eqnarray*}
where
%$\iota : p_i \rightarrow \bP^2$ is the inclusion,
$\scI_{p_i}$ is the ideal sheaf of $p_i$
and the last line is the fiber at $p_i$.
Here, we have used the exact sequence
$$
0 \rightarrow \scO_\dP \rightarrow \scO_\dP(E_i)
  \rightarrow \scO_{E_i}(-1) \rightarrow 0.
$$
We can further use the exact sequence
$$
0 \rightarrow \Omega_{\bP^2}(1) \rightarrow V^\vee \otimes \scO_{\bP^2}
 \rightarrow \scO_{\bP^2}(1) \rightarrow 0
$$
to reduce the computations involving $\Omega_{\bP^2}(1)$
to those involving $\scO_{\bP^2}(i)$, $i \in \bZ$:
\begin{eqnarray*}
\RHom(\scO_{E_i}(-1)[-1], \ectwo)
 & = & \RGamma( \scO_{p_i} \otimes
         \{V^\vee \otimes \scO_{\bP^2}
                \rightarrow \scO_{\bP^2}(1) \} ) \\
 & = & \{V^\vee \rightarrow (\bC p_i)^\vee \},
\end{eqnarray*}
\begin{eqnarray*}
\RHom(\econe, \ectwo)
 & = & \RHom( \scO_{\bP^2}(-1),
         \{ V^\vee \otimes \scO_{\bP^2}
                 \rightarrow \scO_{\bP^2}(1) \} ) \\
 & = & \RGamma( \scO_{\bP^2}(1) \otimes
         \{ V^\vee \otimes \scO_{\bP^2}
                 \rightarrow \scO_{\bP^2}(1) \} ) \\
 & = & \RGamma( \{ V^\vee \otimes \scO_{\bP^2}(1)
                 \rightarrow \scO_{\bP^2}(2) \} ) \\
 & = & \{V^\vee \otimes V^\vee \rightarrow \Sym^2 V^\vee \} \\
 & = & \wedge^2 V^\vee,
\end{eqnarray*}
\begin{eqnarray*}
\RHom(\ectwo, \ecthree)
 & = & \RHom( \{ V^\vee \otimes \scO_{\bP^2}
                 \rightarrow \scO_{\bP^2}(1) \},
               \scO_{\bP^2} ) \\
 & = & \RGamma( \{ \scO_{\bP^2}(-1) \rightarrow
                  V \otimes \scO_{\bP^2} \} ) \\
 & = & V.
\end{eqnarray*}
Compositions of morphisms can be
easily read off
from the above computations.
%Compositions involving sheaves
%supported on an exceptional divisor
%is just the restriction of global Hom's
%between vector bundles on $\bP^2$
%to the fiber above the center of the blow-up:
%The composition
%$$
% \Hom(\scO_{E_i}(-1), \Lphi \scE) \times
% \Hom(\Lphi \scE, \Lphi \scF)
%  \rightarrow \Hom(\scO_{E_i}(-1), \Lphi \scF)
%$$
%gives, for any $\varphi \in \Hom(\scE, \scF)
%\cong \Hom(\Lphi \scE, \Lphi \scF)$,
%a homomorphism
%$$
%\Hom(\scO_{E_i}(-1), \Lphi \scE)
%  = \scE |_{p_i}
% \rightarrow \scF |_{p_i}
%  = \Hom(\scO_{E_i}(-1), \Lphi \scF)
%$$
%which is just the restriction of
%$\Hom(\scE, \scF)$
%to the fiber over $p_i$, and
%the composition
%$$
%\Hom(\econe,\ectwo) \times \Hom(\ectwo, \ecthree)
% \rightarrow \Hom(\econe, \ecthree)
%$$
%is given by the interior product
%$\Lambda^2 V^\vee \times V \rightarrow V^\vee.$
\qed

\section{Fukaya category} \label{Fuk}

First we recall the definition of an $A_\infty$-category.
For a $\bZ$-graded vector space $V$
and $i \in \bZ$,
$V[i]$ denotes the shift of $V$ by $i$;
$V[i]^j=V^{i+j}.$

\begin{definition} \label{def:a_infinity}
An $A_\infty$-category $\scA$ consists of
\begin{enumerate}
\item 
the set of objects $\Ob(\scA)$,
\item
for any $c_1,\; c_2 \in \Ob(\scA)$,
a $\bZ$-graded $\bC$-vector space $\hom_{\scA}(c_1, c_2)$
called
the set of morphisms,
\item
for any positive integer $k$
and for any set of objects $\{c_i\}_{i=0}^k$,
the composition
$$\
 \m_k : \hom_{\scA}(c_0,c_1)[1]\otimes \cdots
                \otimes \hom_{\scA}(c_{k-1},c_k)[1]
                   \longrightarrow \hom_{\scA}(c_0,c_k)[1]
$$
which is a linear map
of degree $1$
\end{enumerate}
such that
for any positive integer $k$,
any set of objects $\{c_i\}_{i=0}^k$
and any set of morphisms $\{a_i\}_{i=1}^k$,
$a_i \in \hom_{\scA}(c_{i-1},c_i)$,
the following {\em $A_\infty$-relations} holds:
\begin{eqnarray*} \label{eq:a_infinity}
 \sum_{i=0}^{k-1} \sum_{j=i+1}^k
  (-1)^{\deg a_1 + \cdots + \deg a_i}
  \m_{i-j+k+1}(a_1 \otimes \cdots \otimes a_i \\
   \otimes 
    \m_{j-i}(a_{i+1} \otimes \cdots \otimes a_j)
   \otimes
    a_{j+1} \otimes \cdots \otimes a_k ) = 0.
\end{eqnarray*}
Here degrees are counted after shifts, i.e., if $a \in V^i$,
$\deg a = i-1$ in $V[1]$.
Since $\m_1^2=0$, we define
$$
 \Hom_{\scA}(c_1,c_2) = H^0(\hom_{\scA}(c_1,c_2),\m_1).
$$
\end{definition}

%\begin{definition}
%Let $\scA$ be an $A_\infty$-category
%and $(c_1, \ldots, c_N)$ be a finite ordered subset
%of $\Ob(\scA)$ such that $\hom_{\scA}(c_i,c_j)$ is
%finite-dimensional for $i<j$.
%The {\em directed $A_\infty$-category} $\scA^{\rightarrow}$
%associated to $(c_1, \ldots, c_N)$
%is an $A_\infty$-category such that
%$$
% \Ob(\scA^{\rightarrow}) = (c_1, \ldots, c_N)
%$$
%and whose morphisms satisfy
%$$
% \hom_{\scA^{\rightarrow}}(c_i,c_j) = \left\{
%  \begin{array}{lc}
%   0 & i>j, \\
%   \bC \cdot \id_{c_i} & i=j, \\
%   \hom_{\scA}(c_i,c_j) & i<j
%  \end{array}
%  \right.
%$$
%and compositions satisfy
%$$
%\m^{\scA^\rightarrow}_k(a_1,\ldots,a_k) = \m^\scA_k(a_1,\ldots,a_k)
%$$
%where $a_j \in \hom_\scA(c_{i_j-1},c_{i_j})$ for $i_o<\cdots<i_k$
%and $\id_{c_i}$ is the strict unit,
%i.e., $\m_k$ involving $\id_{c_i}$ for $k \neq 2$ is zero
%and $\id_{c_i}$ is the unit of $\m_2$.
%\end{definition}

The Fukaya category
of Lagrangian submanifolds
in a symplectic manifold
is defined in
\cite{Fukaya-Oh-Ohta-Ono}.
We use the following adaptation
for exact Morse fibrations
by Seidel \cite{Seidel_VC}.
Let $W: Z \rightarrow \bC$
be a regular function
on an affine algebraic manifold $Z$
of complex dimension $n$
with a K\"{a}hler structure.
%We specialize to the case
%of $\dim Z=2$
%to simplify the exposition.
Assume the following conditions:
\begin{itemize}
 \item The K\"{a}hler metric is complete.
 \item The symplectic form $\omega$
  of $Z$ is exact,
  i.e., there exists a one form $\theta$ on $Z$
  such that $\omega = d \theta$.
 \item At any critical points of $W$,
  the Hessian of $W$ is non-degenerate.
 \item All the critical values are distinct.
\end{itemize}
Such $W$ gives rise to an exact Morse fibration
in the terminology of \cite{Seidel_VC}.
Using the K\"{a}hler structure,
we can define the lift
$\widetilde{c}_p : [0,1] \rightarrow Z$
of a path 
$c : [0,1] \rightarrow \bC$
starting from a point $p \in Z$
such that $W(p) = c(0)$
by using the horizontal distribution
defined as the orthogonal complement
of the tangent space
along the fiber of $W$.

Assume that the origin is
a regular value of $W$ and
fix an order $(p_i)_{i=1}^N$
on the set of critical points of $W$.
A distinguished set $(c_i)_{i=1}^N$
of vanishing paths is a set of smooth paths
$c_i : [0,1] \rightarrow \bC$
satisfying
 \begin{enumerate}
  \item $c_i(0) = 0$, $c_i(1) = W(p_i)$,
  \item $c_i$ has no self-intersection,
  \item images of $c_i$ and $c_j$ intersects only at the origin,
  \item $c_i'(0) \neq 0$, and
  \item $\arg c_{i+1}'(0) < \arg c_i'(0)$, $i=1,\ldots N-1$,
         for some choice of the branch of $\arg(\bullet)$.
 \end{enumerate}
Given a distinguished set
$(c_i)_{i=1}^N$ of vanishing paths,
%from zero to the critical values of $W$,
the corresponding vanishing cycles
$(C_i)_{i=1}^N$
are defined by
$$
 C_i = \{p \in W^{-1}(0) \suchthat
           \lim_{t \rightarrow 1} \widetilde{c}_p(t) = p_i\}.
$$
They are Lagrangian submanifolds
of $W^{-1}(0)$.

The {\em directed Fukaya category} $\FukW$
for $W$ is roughly an $A_\infty$-category
whose objects are vanishing cycles
and whose morphisms are Lagrangian intersection
Floer complexes.
To define a $\bZ$-grading
on the Floer complex,
we need the concept of {\em grading}
on Lagrangian submanifolds
introduced by Kontsevich \cite{Kontsevich_HAMS},
which we now recall.
See also Seidel \cite{S}.

Let $(M,\omega)$ be a symplectic manifold
of dimension $2n$.
An almost complex structure on $M$
is a section $J \in \Gamma(M, \End(T M))$
such that $J^2 = -\id$.
$J$ is called $\omega$-compatible
if $g(V_1,V_2) = \omega(V_1, J V_2)$
%for two vector fields $V_1$ and $V_2$
defines a Hermitian metric on the tangent bundle.
Fix an $\omega$-compatible almost complex structure $J$
of $M$ and
let $\scS$ be the principal $U(1)$-bundle
associated to the complex line bundle
$(\Lambda^n (T^*M, J))^{\otimes 2}$.
The fiber of $\scS$ at $p \in M$ is
$(\Lambda^n (T_p^* M, J))^{\otimes 2} / \bR^{> 0}$.
We assume that the first Chern class
of this complex line bundle vanishes,
so that $\scS$ has a section.
A grading of $M$ is
a choice of a section $\Theta : M \rightarrow \scS$.
Fix a grading $\Theta$ on $M$.
Let $\Lag_M \rightarrow M$ be
the Lagrangian Grassmannian bundle on $M$,
whose fiber at $p \in M$
is the Grassmannian of Lagrangian subspaces
in the symplectic vector space $T_p M$.
Define $\det_\Theta^2 : \Lag_M \rightarrow U(1)$ as follows:
For a Lagrangian subspace $L \subset T_p M$,
pick any basis $\{e_i\}_{i=1}^n$ of $L$
and take the square of their exterior product;
$(e_1 \wedge e_2 \wedge \cdots \wedge e_n)^{\otimes 2}
\in (\Lambda^n (T_p M, J))^{\otimes 2}$.
$\det_\Theta^2(L)$ is the image of this element by $\Theta(p)$.
A Lagrangian submanifold $L \subset M$
gives a canonical section
$s_L : L \ni p \mapsto T_p L \in \Lag_M|_p$.
Denote the composition of $s_L$ and $\det_\Theta^2$ by $\phi_L$.
A grading of a Lagrangian submanifold
is a lift $\widetilde{\phi_L} : L \rightarrow \bR$
of $\phi_L$ to the universal cover $\bR$ of
$U(1) \cong \bR/\bZ$.

Now we define the Maslov index.
A smooth path
$\Lambda : [0,1] \rightarrow \Lag(V,\omega)$
in the Lagrangian Grassmannian
of a fixed symplectic vector space $(V, \omega)$
of dimension $2 n$
is called {\em crossingless}
if $\Lambda(0) \cap \Lambda(t) = \Lambda(0) \cap \Lambda(1)$
for all $t \in (0,1]$.
For a crossingless path $\Lambda_t$,
its differential $\Lambda'(0)$ at $t=0$
gives an element of the tangent space
$T_{\Lambda(0)} \Lag(V,\omega) \subset T_{\Lambda(0)} \Gr(n,V)
  \cong \Hom(\Lambda(0),V/\Lambda(0))$
where $\Gr(n,V)$ is the Grassmannian of $n$-dimensional
subspaces in $V$.
The composition of $\Lambda'(0)$ and $\omega$
defines a quadratic form
$\Lambda(0) \ni v \mapsto \omega(v, \Lambda'(0) v)$
on $\Lambda(0)$,
which descends to a quadratic form on 
$\Lambda(0)/\Lambda(0) \cap \Lambda(1)$
since it vanishes on $\Lambda(0) \cap \Lambda(1)$.
The resulting form on $\Lambda(0)/\Lambda(0) \cap \Lambda(1)$
is called the {\em crossing form} \cite{Robbin-Salamon}.
For an intersection $p \in L_1 \cup L_2$
of two graded Lagrangian submanifolds
$(L_1, \widetilde{\phi_{L_1}})$ and $(L_2, \widetilde{\phi_{L_2}})$,
its Maslov index is defined as follows:
Choose a crossingless path
$\Lambda : [0,1] \rightarrow \Lag_M|_p$
from $T_p L_0$ to $T_p L_1$
such that the corresponding crossing form
at $t=0$ is negative definite.
There is a unique lift 
$\widetilde{\alpha} : [0,1] \rightarrow \bR$
of the composition $\alpha$ of
$\Lambda : [0,1] \rightarrow \Lag_M|_p$ and
$\det_\Theta^2(p) : \Lag_M|_p \rightarrow U(1)$
such that $\widetilde{\alpha}(0)=\widetilde{\phi_{L_0}}(p)$,
and the Maslov index $I(p)$ of the intersection
$p$ of two graded Lagrangian submanifolds
is defined by
$$
 I(p) = \widetilde{\phi_{L_1}}(p) - \widetilde{\alpha}(1).
$$

Now we come back to our exact Morse fibration $W$.
A relative Maslov map
is a section $\Theta$
of the second tensor power
$(\Omega_{Z / \bC}^{n-1})^{\otimes 2}$
of the top exterior product
of the relative cotangent bundle
away from the critical points.
Since $0$ is a regular value,
the restriction of $\Theta$
to $W^{-1}(0)$ gives a grading
of $W^{-1}(0)$.
Assume that
all the vanishing cycles can be graded,
and fix a grading on each vanishing cycle.
We also need to choose a spin structure
on each vanishing cycles
in order to orient
the moduli spaces of
pseudoholomorphic maps
\cite{Fukaya-Oh-Ohta-Ono}.
Since
each vanishing cycle is homeomorphic
to a sphere,
the choice of a spin structure
is unique
except for $\dim_\bC W^{-1}(0) = 1$,
where one has as many as
$H^1(S^1,\bZ/2 \bZ)=\bZ / 2 \bZ$
choices.
We take the non-trivial spin structure
in such a case.
Let $C_i^\flat$ denote
the vanishing cycle $C_i$
endowed with the above grading and spin structure.
We assume that vanishing cycles intersects each other
transversally.
This condition can always be met
by moving vanishing cycles
within their Hamiltonian isotopy classes
if necessary.

\begin{definition}
Given a function $W$ on an affine K\"{a}hler manifold $Z$
together with a relative Maslov map $\Theta$
and a choice of a distinguished basis
of vanishing cycles
$(C_i^\flat)_{i=1}^N$
with gradings and spin structures,
its directed Fukaya category $\FukW$
is an $A_\infty$-category such that
\begin{itemize}
 \item the set of objects
        is the distinguished basis of vanishing cycles;
        $$ \Ob(\FukW) = (C_1^\flat, \ldots, C_N^\flat),$$
 \item the set of morphisms between
        $C_i^\flat$ and $C_j^\flat$ is
        the $\bZ$-graded vector space
$$
\hFukW(C_i^\flat, C_j^\flat)
 = \left\{
  \begin{array}{lc}
   0 & i>j, \\
   \bC \cdot \id_{C_i^\flat} & i=j, \\
   \bigoplus_{p \in C_i \cap C_j} \vspan_{\bC}\{p\} & i<j
  \end{array}
  \right.
$$
      where $\deg p= I(p)$ (the Maslov index), and
 \item for a positive integer $k$,
        a set of objects
        $(C_{i_0}^\flat,\ldots,C_{i_k}^\flat)$
        and morphisms $p_l \in C_{i_{l-1}} \cap C_{i_l}$
        for $l = 1,\ldots,k$,
        the composition $\m_k$ is given by
$$
 \m_k(p_1,\ldots,p_k) = \sum_{p_0 \in C_{i_0} \cap C_{i_k}}
     \# \M_{k+1}(C_{i_0},\ldots,C_{i_k};p_0,\ldots,p_k) p_0.
$$
\end{itemize}
\end{definition}
Here, $\M_{k+1}(C_{i_0},\ldots,C_{i_k};p_0,\ldots,p_k)$
is the stable compactification of \\
$\scM_{k+1}(C_{i_0},\ldots,C_{i_k};p_0,\ldots,p_k)$
defined as follows:
A disk with $k+1$ marked points
on the boundary
is a pair $(D^2,(z_0,\ldots,z_k))$
of a closed unit disk
$D^2=\{z \in \bC \suchthat |z| \leq 1\}$
and an ordered set
$(z_0,\ldots,z_k)$
of $k+1$ points on the boundary
respecting the cyclic order.
Let $\partial_l D^2 \in \partial D^2$
be the interval
between $z_l$ and $z_{l+1}$,
where we set $z_{k+1} = z_0$.
Fix an almost complex structure $J$ on $W^{-1}(0)$.
A smooth map $\varphi : D^2 \rightarrow M$
is called pseudoholomorphic if
$$
d \varphi \circ J_{D^2} = J \circ d \varphi
$$
where $J_{D^2}$ is the canonical complex structure
on $D^2$.
$\scM_{k+1}(C_{i_0},\ldots,C_{i_k};p_0,\ldots,p_k)$ is
the moduli space of pairs
$((D^2,(z_0,\ldots,z_k)),\varphi)$
such that
\begin{enumerate}
 \item $(D^2,(z_0,\ldots,z_k))$ is a disk
   with $k+1$ marked points on the boundary,
 \item $\varphi : D^2 \rightarrow W^{-1}(0)$
   is a pseudoholomorphic map,
 \item $\varphi(\partial_l D^2) \subset C_{i_l}$ for $l=0,\ldots,k$, and
 \item $\varphi(z_l) = p_l$ for $l=0,\ldots,k$.
\end{enumerate}

Although $\M_{k+1}(C_{i_0},\ldots,C_{i_k};p_0,\ldots,p_k)$
is not an honest manifold in general
due to subtleties in the definition of
the moduli space of pseudoholomorphic maps,
it has a Kuranishi structure with corners
and it makes sense to ``count'' numbers
$\# \M_{k+1}(C_{i_0},\ldots,C_{i_k};p_0,\ldots,p_k)$
of points.
See \cite{Fukaya-Oh-Ohta-Ono}
for details.
These numbers are counted with signs
determined by the orientation of the moduli space,
which in turn is determined
by the spin structures on vanishing cycles.
$\# \M_{k+1}(C_{i_0},\ldots,C_{i_k};p_0,\ldots,p_k)$
is zero if \\
$\dim \M_{k+1}(C_{i_0},\ldots,C_{i_k};p_0,\ldots,p_k) \neq 0$.

These moduli spaces must be compact
so that it makes sense
to count the numbers of points.
To prove compactness of the moduli space,
one needs the boundedness
of the energy.
Usually, this is achieved
by introducing
the Novikov ring
as the coefficient ring
of the Floer complex
in order to control energies
of pseudoholomorphic maps.
This is not necessary in our case
since vanishing cycles are exact Lagrangian submanifolds,
i.e., $[\theta] = 0 \in H^1(C_i,\bR)$
for the primitive $\theta$ of $\omega$;
$\omega = d \theta$.
To see this,
let $C$ be the vanishing cycle
corresponding to a vanishing path
$c : [0,1] \rightarrow \bC$.
Then $C$ bounds the disk
$D = \bigcup_{p \in C} \widetilde{c_p}([0,1])$;
$\partial D = C$.
Since
$d \theta|_D = \omega|_D = 0$,
$[\theta] = 0 \in H^1(\partial D,\bR)$.
Therefore,
there exists a function $K_i$
on each vanishing cycles
such that
$\theta|_{C_i} = d K_i$.
Then for any $\varphi \in
\M_{k+1}(C_{i_0},\ldots,C_{i_k};p_0,\ldots,p_k)$,
the energy of $\varphi$ is
\begin{eqnarray*}
 E(\varphi) & = & \int_{D^2} \varphi^*\omega
  =  \int_{\partial D^2} \varphi^* \theta
  =  \sum_{l=0}^k \int_{\partial_i D^2} \varphi^*\theta \\
 &=& \sum_{l=0}^k (K_{i_l}(p_{l+1})- K_{i_l}(p_l))
\end{eqnarray*}
and does not depend on the homotopy class of $\varphi$.
This also implies the vanishing of $\m_0$.
Non-compactness of $W^{-1}(0)$
does not cause any problem either.
Since $W^{-1}(0)$ is Stein,
the maximum principle prevents
pseudoholomorphic  disks
from running away to infinity,
assuring the compactness of the moduli space.

Seidel proved that although
the directed Fukaya category $\FukW$
depends on the choice of the distinguished basis
of vanishing cycles,
different choices are related by mutations,
hence its derived category
is invariant.

Now we explain the mirror construction
of toric Fano manifolds
after Givental \cite{Givental_HGMS}.
Given a fan of an $n$-dimensional toric Fano manifold,
let $\{v_i\}_{i=1}^r$,
$v_i = (v_{i1},\ldots,v_{in}) \in \bZ^n$,
be the set of generators of
its one-dimensional cones.
Then the mirror object
for this
toric Fano manifold is
the regular function
\begin{equation} \label{eq:generalW}
 W(x_1,\ldots,x_n)=\sum_{i=1}^r q_i x_1^{v_{i1}} \cdots x_n^{v_{in}}
\end{equation}
on the algebraic torus $\Spec \bC[x_i^{\pm 1}]_{i=1}^n$.
Here, $q_i$'s are parameters
corresponding to the deformation
of symplectic structures on the toric Fano manifold.

\begin{figure}[htbp]
\centering
\psfrag{c1}{$c_1$}
\psfrag{c2}{$c_2$}
\psfrag{c3}{$c_3$}
\psfrag{c4}{$c_4$}
\psfrag{c5}{$c_5$}
\psfrag{c6}{$c_6$}
\includegraphics{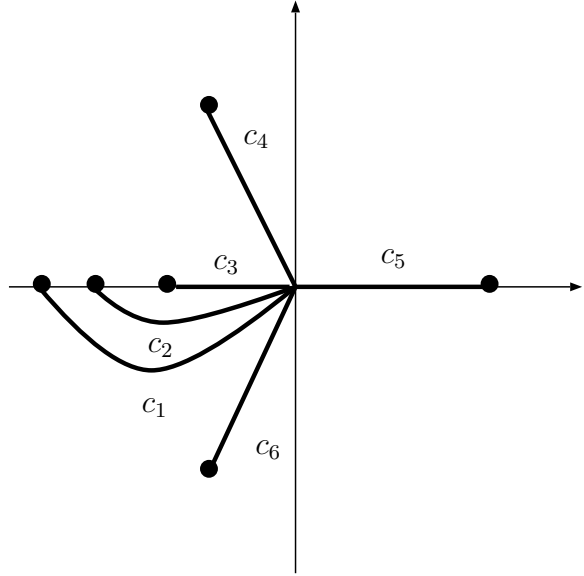}
\caption{Distinguished basis of vanishing paths}
\label{fg:vanishing_paths}
\end{figure}

\begin{figure}[htbp]
\centering
\psfrag{C3}{$C_1$}
\psfrag{C4}{$C_2$}
\psfrag{C5}{$C_3$}
\psfrag{C6}{$C_4$}
\psfrag{C7}{$C_5$}
\psfrag{C8}{$C_6$}
\psfrag{x1}{$x_1$}
\psfrag{x2}{$x_2$}
\psfrag{x3}{$x_3$}
\psfrag{x12*}{$x_1^\vee \wedge x_2^\vee$}
\psfrag{x23*}{$x_2^\vee \wedge x_3^\vee$}
\psfrag{x31*}{$x_3^\vee \wedge x_1^\vee$}
\psfrag{x1*}{$x_1^\vee$}
\psfrag{-x1*}{$-x_1^\vee$}
\psfrag{x2*}{$x_2^\vee$}
\psfrag{x3*}{$x_3^\vee$}
\psfrag{e1}{$e_{1}$}
\psfrag{e2}{$e_{2}$}
\psfrag{e3}{$e_{3}$}
\scalebox{1}{\includegraphics{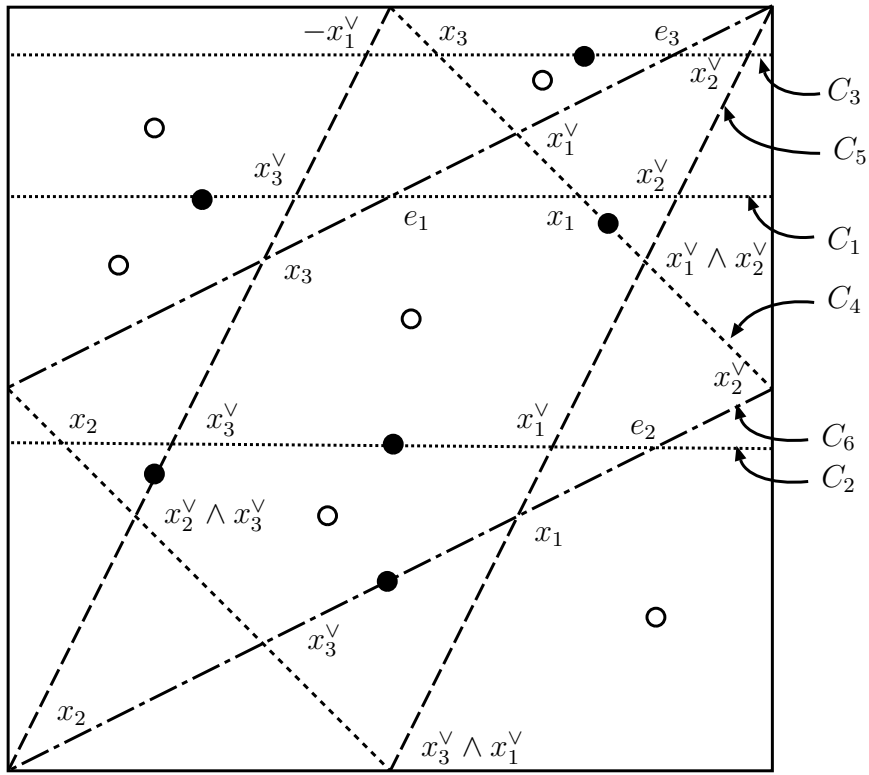}}
\caption{Vanishing cycles in the fiber at the origin}
\label{fg:vanishing_cycles}
\end{figure}

Therefore, the mirror of our $\dP$
is the regular function
\begin{equation} \label{eq:W}
 W(x,y) = q_1 x + q_2 y + \frac{q_3}{x y}
          + \frac{q_4}{x} + \frac{q_5}{y} + q_6 x y 
\end{equation}
on the algebraic torus 
$(\bC^\times)^2 = \Spec \bC[x,x^{-1},y,y^{-1}]$
(See Figure \ref{fg:dp3}).
The Fukaya category does not depend
on a general choice of $q_i$'s.
We have used
$(q_1,q_2,q_3,q_4,q_5,q_6) = (1,1,1,0.215,0.25,0.3)$
to draw the figures appearing below.
Equip $(\bC^\times)^2$
with the symplectic form
$\frac{d|x|}{|x|} \wedge d(\arg x) + \frac{d|y|}{|y|} \wedge d(\arg y)$
and
the relative Maslov map
$$
 \Theta = (\Res \frac{dx \wedge dy}{x y W(x,y)})^{\otimes 2}
        = (\frac{dx}{\partial_y (x y W(x,y))})^{\otimes 2}.
$$
$W^{-1}(0)$ is an affine elliptic curve,
which can be compactified by adding six points.
We depict the critical values of $W$
and our choice of a distinguished set
of vanishing paths
in Figure \ref{fg:vanishing_paths}.
Figure \ref{fg:vanishing_cycles}
shows the corresponding
%distinguished basis of
vanishing cycles.

The opposite sides of the square
in Figure \ref{fg:vanishing_cycles}
are identified
to form a two-dimensional torus.
The open circles denote the points
which are missing
due to the non-compactness
of $(\bC^\times)^2$.
We have drawn Figure \ref{fg:vanishing_cycles}
in the following way:
First regard $W^{-1}(t)$ as a branched double
cover of $\bC^\times$
by the projection
$\pi_t : W^{-1}(t) \ni (x,y) \mapsto y \in \bC^\times$.
For any $i =1,\ldots,6$,
the branch points of $\pi_{c_i(t)}$ moves in $\bC^\times$
as we vary $t$,
until two of them finally collides
at $t=1$.
The vanishing cycle $C_i$ is
the circle in $W^{-1}(0)$
over this trajectory of collision.
This determines the vanishing cycle $C_i$
up to isotopy.
These vanishing cycles
can be straightened within their Hamiltonian isotopy classes.
Here, being straight refers to
the flat metric on the torus
$\overline{W^{-1}(0)}$,
where $\overline{\bullet}$ denotes
the completion.
Note that this flat metric on
$\overline{W^{-1}(0)}$
has nothing to do
with the K\"{a}hler metric
on $W^{-1}(0)$
induced from that of
$(\bC^\times)^2$.
This determines $C_i$ up to translation.
This translational ambiguity can be fixed
by imposing exactness,
which requires the knowledge of the primitive
$\theta$ of the symplectic form $\omega$
on $W^{-1}(0)$.
We do not try to carry this out
since translations of straight $C_i$'s
do not alter the combinatorial structure
of intersections,
which is all we need in our computation of the
Fukaya category below.
The grading on $W^{-1}(0)$
given by $\Theta$
is the grading
coming from the restriction
of the second tensor power
of the holomorphic 1-form
on $\overline{W^{-1}(0)}$.
Then
one can choose a grading $\widetilde{\phi_{C_i}}$
on each vanishing cycle $C_i$
so that
all the Maslov indices
of $C_i \cap C_j$
is zero for $i < j$.
To give a spin structure on $C_i$
is the same as to give a two-fold covering of $C_i$.
We take the non-trivial cover
such that two branches interchange
on the black dots in Figure \ref{fg:vanishing_cycles}.
Different choices for spin structures
lead to different categories,
and it turns out that
the above choice gives a category
derived equivalent to
the category of coherent sheaves
on the del Pezzo surface $\dP$.
Let $C_i^\flat$ denote the vanishing cycle $C_i$
with the above grading and spin structure.
%Different choices for branch points
%are related simply by
%the sign flips
%of the generators
%of the homomorphisms
%corresponding to
%intersection points
%of vanishing cycles.
Since the Maslov index
of all the intersection points are zero,
the sign for the counting of a triangle
is $-1$ if the boundary of the triangle
hits odd numbers of these dots,
and $+1$ otherwise
(see Seidel \cite{Seidel_K3}).
Note that although the choice of spin structures
is important,
the choice of the positions
of the black dots is irrelevant.
The change of the signs
caused by a change of the positions
of the black dots
can be absorbed
by a redefinition of
the signs of the basis for the Floer cohomologies.
$\m_k$ is non-zero only for $k=2$
(polygons with more than three edges
in Figure \ref{fg:vanishing_cycles}
do not contribute
since $\hFukW(C^\flat_i,C^\flat_j) = 0$ for $i > j$).
%By the assignment of objects
%\begin{eqnarray*}
% C_1 & \longleftrightarrow & \scO_{E_1}(-1), \\
% C_2 & \longleftrightarrow & \scO_{E_2}(-1), \\
% C_3 & \longleftrightarrow & \scO_{E_3}(-1), \\
% C_4 & \longleftrightarrow & \econe, \\
% C_5 & \longleftrightarrow & \ectwo, \\
% C_6 & \longleftrightarrow & \ecthree,
%\end{eqnarray*}
%and that of morphisms
%as in Figure \ref{fg:vanishing_cycles},
%one can check that
%the compositions perfectly match
%the table in the Appendix.
%Since both $\DbdP$ and $D^b \Fuk W$
%are generated by exceptional collections
%with identical morphisms and compositions,
%the equivalence in
%Theorem \ref{th:hms} follows.

\begin{figure}[H]
\centering
\psfrag{C3}{$C_1$}
\psfrag{C4}{$C_2$}
\psfrag{C5}{$C_3$}
\psfrag{C6}{$C_4$}
\psfrag{C7}{$C_5$}
\psfrag{C8}{$C_6$}
\psfrag{x1}{$x_1$}
\psfrag{x2}{$x_2$}
\psfrag{x3}{$x_3$}
\psfrag{x12*}{$x_1^\vee \wedge x_2^\vee$}
\psfrag{x23*}{$x_2^\vee \wedge x_3^\vee$}
\psfrag{x31*}{$x_3^\vee \wedge x_1^\vee$}
\psfrag{x1*}{$x_1^\vee$}
\psfrag{x2*}{$x_2^\vee$}
\psfrag{x3*}{$x_3^\vee$}
\psfrag{e1}{$e_{1}$}
\psfrag{e2}{$e_{2}$}
\psfrag{e3}{$e_{3}$}
\scalebox{0.7}{\includegraphics{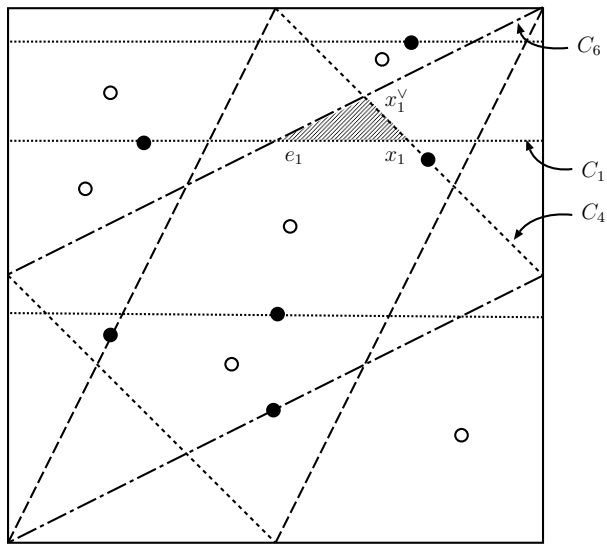}}
\caption{Counting triangle}
\label{fg:cd1}
\end{figure}

\begin{theorem} \label{th:commutativity}
There exists an isomorphism
$$ 
 \phi_{ij} : \Hom_\DbdP(\scE_i, \scE_j)
       \rightarrow \HFukW(C_i^\flat,C_j^\flat)
$$
as a $\bC$-vector space
for $i,j=1,\ldots,6$
such that
the diagrams
\begin{equation} \label{eq:diagram}
\begin{CD}
 \Hom_\DbdP(\scE_i, \scE_j) \times \Hom_\DbdP(\scE_j, \scE_k)
   @>>> \Hom_\DbdP(\scE_i, \scE_k) \\
 @VV{\phi_{ij} \times \phi_{jk}}V @VV{\phi_{ik}}V \\
 \HFukW(C_i^\flat,C_j^\flat) \times \HFukW(C_j^\flat,C_k^\flat)
   @>>> \HFukW(C_i^\flat,C_k^\flat)
\end{CD}
\end{equation}
commute.
\end{theorem}
\proof
We omit the lower suffix of $\Hom(\bullet,\bullet)$
since there is no danger of confusion.
We construct the above $\phi_{ij}$'s explicitly.
Since $\m_1=0$,
$\Hom(C_i^\flat,C_j^\flat)$
is isomorphic to $\hom(C_i^\flat,C_j^\flat)$
and spanned by intersection points
of $C_i$ and $C_j$, to which we assign
the basis of $\Hom(\scE_i,\scE_j)$
as in Figure \ref{fg:vanishing_cycles}.
Notations for the basis of
$\Hom(\scE_i, \scE_j)$
is given in Appendix.
This defines the linear isomorphisms $\phi_{ij}$'s.
The commutativity
of (\ref{eq:diagram})
is verified
by counting triangles.
Let us illustrate this
with a few examples.
Take $x_1 \in \Hom(\scE_1, \scE_4)$
and $x_1^\vee \in \Hom(\scE_4, \scE_6)$,
whose composition is
$e_1 \in \Hom(\scE_1, \scE_6)$.
On the Fukaya side,
we count the number of triangles
whose edges are contained in
$C_1 \cup C_4 \cup C_6$
and two of whose vertices are
$x_1$ and $x_1^\vee$.
Such a triangle exists uniquely,
whose remaining vertex is $e_1$
(Figure \ref{fg:cd1}).
%This matches the computation
%on the derived category of coherent sheaves
%on the del Pezzo surface.

\begin{figure}[H]
\centering
\psfrag{C3}{$C_1$}
\psfrag{C4}{$C_2$}
\psfrag{C5}{$C_3$}
\psfrag{C6}{$C_4$}
\psfrag{C7}{$C_5$}
\psfrag{C8}{$C_6$}
\psfrag{x1}{$x_1$}
\psfrag{x2}{$x_2$}
\psfrag{x3}{$x_3$}
\psfrag{x12*}{$x_1^\vee \wedge x_2^\vee$}
\psfrag{x23*}{$x_2^\vee \wedge x_3^\vee$}
\psfrag{x31*}{$x_3^\vee \wedge x_1^\vee$}
\psfrag{x1*}{$x_1^\vee$}
\psfrag{x2*}{$x_2^\vee$}
\psfrag{x3*}{$x_3^\vee$}
\psfrag{e1}{$e_{1}$}
\psfrag{e2}{$e_{2}$}
\psfrag{e3}{$e_{3}$}
\scalebox{0.7}{\includegraphics{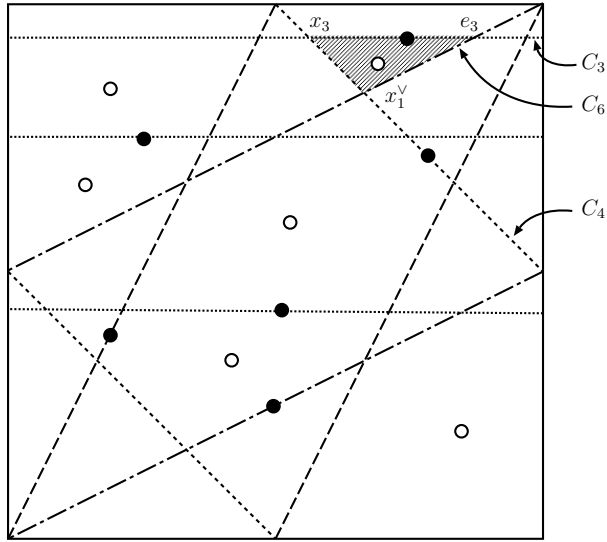}}
\caption{No triangles}
\label{fg:cd2}
\end{figure}

Next we try to compose
$x_3 \in \Hom(C_3^\flat,C_4^\flat)$
and $x_1^\vee \in \Hom(C_4^\flat, C_6^\flat)$
in the Fukaya category.
The only possibility
would be the triangle
with $e_3$ as the remaining vertex.
However, this is not allowed
because of
the missing point
(Figure \ref{fg:cd2}).

\begin{figure}[H]
\centering
\psfrag{C3}{$C_1$}
\psfrag{C4}{$C_2$}
\psfrag{C5}{$C_3$}
\psfrag{C6}{$C_4$}
\psfrag{C7}{$C_5$}
\psfrag{C8}{$C_6$}
\psfrag{x1}{$x_1$}
\psfrag{x2}{$x_2$}
\psfrag{x3}{$x_3$}
\psfrag{x12*}{$x_1^\vee \wedge x_2^\vee$}
\psfrag{x23*}{$x_2^\vee \wedge x_3^\vee$}
\psfrag{x31*}{$x_3^\vee \wedge x_1^\vee$}
\psfrag{x1*}{$x_1^\vee$}
\psfrag{x2*}{$x_2^\vee$}
\psfrag{x3*}{$x_3^\vee$}
\psfrag{e1}{$e_{1}$}
\psfrag{e2}{$e_{2}$}
\psfrag{e3}{$e_{3}$}
\scalebox{0.7}{\includegraphics{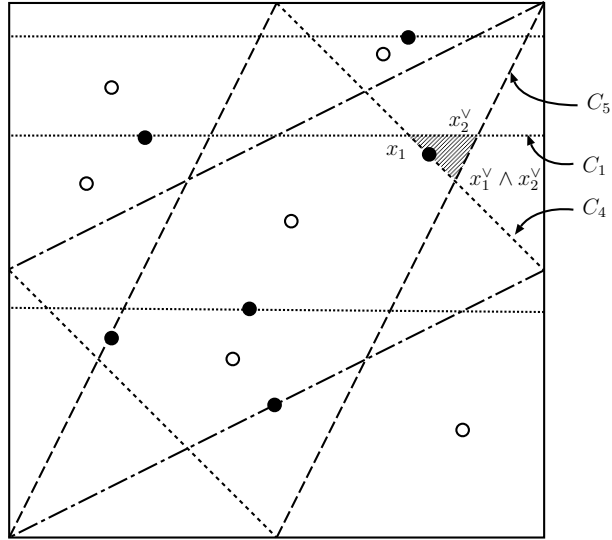}}
\caption{Counting triangles with signs}
\label{fg:cd3}
\end{figure}

This shows that
the composition of $x_3$
and $x_1^\vee$
in the Fukaya category
is zero,
which matches the table
in Appendix
calculated in $\DbdP$.
Finally, as an example of
counting with signs,
we compute
$$
\Hom(C_1^\flat, C_4^\flat) \times
   \Hom(C_4^\flat,C_5^\flat)
 \ni (x_1, x_1^\vee \wedge x_2^\vee)
  \mapsto -x_2^\vee \in \Hom(C_1^\flat, C_5^\flat).
$$
Indeed, the triangle
whose vertices are $x_1, x_1^\vee \wedge x_2^\vee$,
and $x_2^\vee$ contains one black dot 
on its edge
(Figure \ref{fg:cd3}).
We can do similar analyses
for the rest of the triangles.
The result perfectly agrees with the computations
in $\DbdP$.
\qed

Although the derived category of
an $A_\infty$-category
is usually defined
using twisted complexes,
we adopt the following definition
in this paper
since our $\FukW$ satisfies
$\hFukW^k(C_i^\flat,C_j^\flat)=0$
for $k \neq 0$
and $\m_k=0$
for $k \neq 2$.
\begin{definition}
Let $A=\oplus_{i,j}\HFukW(C_i^\flat,C_j^\flat)$
be the total morphism algebra.
The derived Fukaya category
$D^b \FukW$ is
the bounded derived category
$D^b (\modm A)$
of the category of right
finite-dimensional modules
over the algebra $A$.
\end{definition}
Now we can state our main theorem:
\begin{theorem} \label{th:main_theorem}
There exists an equivalence of triangulated categories
$$
\DbdP \cong D^b \FukW.
$$
\end{theorem}
\proof
From Theorem \ref{th:commutativity},
we have
$$
A \cong \oplus_{i,j} \Hom_\DbdP(\scE_i,\scE_j).
$$
Theorem \ref{th:main_theorem} follows
immediately from the theorem of Bondal
\cite{Bondal}
that $\DbdP \cong D^b(\modm (\oplus_{i,j} \Hom(\scE_i,\scE_j)))$.
\qed

We can perform
similar analyses
for all the other toric del Pezzo surfaces
to obtain Theorem \ref{th:hms}.

\section{Appendix}

Here we give the table of compositions of morphisms
in the derived category of coherent sheaves
on the projective plane
$\bP(V)$
blown-up at
$p_1 = [1:0:0]$, $p_2 = [0:1:0]$, $p_3 = [0:0:1]$.
We use the following notations:
$\scE_i$'s are defined in (\ref{eq:ec}),
$\{x_i\}_{i=1}^3$ is the basis of $V$,
$\{x_i^\vee\}_{i=1}^3$ is the dual basis of $V^\vee$, and
$e_i$ is the generator of
$\Hom(\scE_i, \scE_6)=\bC$ for $i=1,2,3$.

$\Hom(\eeone, \eefour)
 \times \Hom(\eefour, \eefive)
 \rightarrow \Hom(\eeone, \eefive)$
$$
\begin{array}{c|c|c|c}
 & x_1^\vee \wedge x_2^\vee
 & x_2^\vee \wedge x_3^\vee 
 & x_3^\vee \wedge x_1^\vee \\
 \hline
 x_1 & -x_2^\vee & 0 & x_3^\vee
\end{array}
$$
$\Hom(\eeone, \eefour) \times \Hom(\eefour, \eesix)
 \rightarrow \Hom(\eeone, \eesix)$
$$
\begin{array}{c|c|c|c}
 & x_1^\vee & x_2^\vee & x_3^\vee \\
 \hline
 x_1 & e_1 & 0 & 0
\end{array}
$$
$\Hom(\eeone, \eefive) \times \Hom(\eefive, \eesix)
 \rightarrow \Hom(\eeone, \eesix)$
$$
\begin{array}{c|c|c|c}
 & x_1 & x_2 & x_3 \\
 \hline
 x_2^\vee & 0 & e_1 & 0 \\
 \hline
 x_3^\vee & 0 & 0 & e_1 \\
\end{array}
$$
$\Hom(\eetwo, \eefour) \times \Hom(\eefour, \eefive)
 \rightarrow \Hom(\eetwo, \eefive)$
$$
\begin{array}{c|c|c|c}
 & x_1^\vee \wedge x_2^\vee
 & x_2^\vee \wedge x_3^\vee 
 & x_3^\vee \wedge x_1^\vee \\
 \hline
 x_2 & x_1^\vee & -x_3^\vee & 0
\end{array}
$$
$\Hom(\eetwo, \eefour) \times \Hom(\eefour, \eesix)
 \rightarrow \Hom(\eetwo, \eesix)$
$$
\begin{array}{c|c|c|c}
 & x_1^\vee & x_2^\vee & x_3^\vee \\
 \hline
 x_2 & 0 & e_2 & 0
\end{array}
$$
$\Hom(\eetwo, \eefive) \times \Hom(\eefive, \eesix)
 \rightarrow \Hom(\eetwo, \eesix)$
$$
\begin{array}{c|c|c|c}
 & x_1 & x_2 & x_3 \\
 \hline
 x_1^\vee & e_2 & 0 & 0 \\
 \hline
 x_3^\vee & 0 & 0 & e_2
\end{array}
$$
$\Hom(\eethree, \eefour) \times \Hom(\eefour, \eefive)
 \rightarrow \Hom(\eethree, \eefive)$
$$
\begin{array}{c|c|c|c}
 & x_1^\vee \wedge x_2^\vee
 & x_2^\vee \wedge x_3^\vee 
 & x_3^\vee \wedge x_1^\vee \\
 \hline
 x_3 & 0 & x_2^\vee & -x_1^\vee
\end{array}
$$
$\Hom(\eethree, \eefour) \times \Hom(\eefour, \eesix)
 \rightarrow \Hom(\eethree, \eesix)$
$$
\begin{array}{c|c|c|c}
 & x_1^\vee & x_2^\vee & x_3^\vee \\
 \hline
 x_3 & 0 & 0 & e_3
\end{array}
$$
$\Hom(\eethree, \eefive) \times \Hom(\eefive, \eesix)
 \rightarrow \Hom(\eethree, \eesix)$
$$
\begin{array}{c|c|c|c}
 & x_1 & x_2 & x_3 \\
 \hline
 x_1^\vee & e_2 & 0 & 0 \\
 \hline
 x_2^\vee & 0 & e_3 & 0
\end{array}
$$
$\Hom(\eefour, \eefive) \times \Hom(\eefive, \eesix)
 \rightarrow \Hom(\eefour, \eesix)$
$$
\begin{array}{c|c|c|c}
 & x_1 & x_2 & x_3 \\
 \hline
 x_1^\vee \wedge x_2^\vee & x_2^\vee & -x_1^\vee & 0 \\
 \hline
 x_2^\vee \wedge x_3^\vee & 0 & x_3^\vee & -x_2^\vee \\
 \hline
 x_3^\vee \wedge x_1^\vee & -x_3^\vee & 0 & x_1^\vee \\
\end{array}
$$

\bibliographystyle{plain}
\bibliography{bibs}

Research Institute of Mathematical Sciences,
Kyoto University,
Oiwake-cho,
Kitashirakawa,
Sakyo-ku,
Kyoto,
606-8502,
Japan.

{\em e-mail address}\ : \  kazushi@kurims.kyoto-u.ac.jp

\end{document}